\documentclass{article}

\usepackage{amsmath,amssymb}
\usepackage{graphicx}
\usepackage{color}
\usepackage{dsfont}

\newtheorem{theorem}{Theorem}[section]%
\newtheorem{corollary}[theorem]{Corollary}%
\newtheorem{lemma}[theorem]{Lemma}%
\newtheorem{proposition}[theorem]{Proposition}%

\newtheorem{definition}[theorem]{Definition}%
\newtheorem{remark}[theorem]{Remark}%
\newtheorem{example}[theorem]{Example}%

\newcommand{\N}{\ensuremath{\mathbb N}}

\newcommand{\R}{\ensuremath{\mathbb R}}
\newcommand{\C}{\ensuremath{\mathbb C}}

\newcommand{\done}{\hfill $\Box$ }

\newcommand{\ls}[1]
    {\dimen0=\fontdimen6\the\font\lineskip=#1\dimen0
     \advance\lineskip.5\fontdimen5\the\font
     \advance\lineskip-\dimen0
     \lineskiplimit=0.9\lineskip
     \baselineskip=\lineskip
     \advance\baselineskip\dimen0
     \normallineskip\lineskip\normallineskiplimit\lineskiplimit
     \normalbaselineskip\baselineskip
     \ignorespaces}

\begin{document}

\bibliographystyle{abbrv}

\title{Some results on the Orlicz space generated from a random normed module}
\author{Mingzhi Wu$^1$ \quad Long Long$^{2,*}$\quad Xiaolin Zeng$^3$\\
1. School of Mathematics and Physics, China University of Geosciences,\\ Wuhan {\rm 430074}, China \\
Email: wumz@cug.edu.cn\\
2. School of Mathematics and Statistics,
Central South University,\\
Changsha {\rm 410083}, China \\
Email: longlong@csu.edu.cn\\
3.School of Mathematics and Statistics, Chongqing Technology and Business University, \\
Chongqing {\rm 400067}, China\\
Email: xlinzeng@163.com\\}

\date{}
 \maketitle

\thispagestyle{plain}
\setcounter{page}{1}

\begin{abstract}
 Noting the important role the abstract $L^p$ space has played in the development of random normed modules, in this paper we introduce and study the Orlicz space generated from a random normed module. First, we give a basic dual space representation theorem which identify the dual of the Orlicz heart of a random normed module with the Orlicz space generated from the random conjugate space. Then, we establish the respective equivalence relations of the strict convexity and uniform convexity of this abstract Orlicz space to the random strict convexity and random uniform convexity of the underlying random normed module. These results demonstrate that it is possible to use the Orlicz space theory in the  further development of random nomed modules.\\


\end{abstract}

\ls{1.5}
\section{Introduction}
The notion of a random normed module (briefly, an RN module) is a generalization of that of a normed space. The theory of random conjugate spaces for RN modules is crucial for the deep development of RN modules. In the last 20 plus years, both the theory of RN modules together with their random conjugate spaces and its applications have undergone a systematic and deep
development (a comprehensive list of main results is included in \cite{Guo-JFA,Guo-recent}), in particular RN modules have been proved to be natural and universal model spaces for conditional risk measures in the recent years, see \cite{FKV1,FKV2,Guo-recent,GZZ}.

In the course of the development of RN modules, the abstract $L^p$-space generated from an RN module is a powerful tool. Let $E$ be an RN module and $p\in [1,+\infty]$, then the $L^p$-space generated from $E$, denoted by $L^p(E)$, consists of the element in $E$ whose random norm is $L^p$-integrable (or essentially bounded if $p=+\infty$). Under the apparent norm induced from $L^p$, $L^p(E)$ is a normed space. Guo established in \cite{Guo-Xiamen,Guo-rep} a basic representation theorem which says that the classical conjugate space of $L^p(E)$ is isometrically and isomorphism to $L^q(E^\ast)$, namely the $L^q$-space generated from $E$'s random conjugate space $E^\ast$, where $1\leq p<+\infty$ and $q$ is the H\"{o}lder conjugate number of $p$.
This theorem unifies all the dual representation theorems
of Lebesgue-Bochner function spaces (see \cite{Guo-rep}). For the development of RN modules, this theorem bridges the $L^p$-space theory with the study of RN modules. Making use of
this theorem, Guo and Li \cite{GL} proved the James theorem in complete RN modules; Guo, Xiao and Chen \cite{GXC} established
a basic strict separation theorem in random locally convex modules; Zhang and Guo \cite{ZG} established a mean ergodic
theorem on random reflexive RN  modules, and Wu \cite{Wu} proved the Bishop-Phelps theorem in complete RN modules endowed with the $(\varepsilon, \lambda)$-topology.

It is well-known that the notion of Orlicz spaces is a generalization of $L^p$ spaces. Orlicz spaces share many useful properties with $L^p$ spaces, among which the most important is that, they are Banach spaces and admit nice duality. Please refer to \cite{ES,RR} for the Orlicz space theory. Considering the applications mentioned above of the abstract $L^p$-space to the study of RN modules, we naturally introduce and study the Orlicz space generated from an RN module.

In this paper, we first introduce the notion of the Orlicz space generated from an RN module. Then, we give a basic dual space representation theorem which identify the dual of the Orlicz heart of an RN module with the Orlicz space generated from the random conjugate space. Finally, we establish the respective equivalence relations of the strict convexity and uniform convexity of this Orlicz space to the random strict convexity and random uniform convexity of the underlying RN module. These results demonstrate that it is possible to use the Orlicz space theory to the study of RN modules in the future.

\section{Terminology and notation}

Let $(\Omega,{\mathcal F},P)$ be a probability space, $K$ the scalar field of real numbers $\R$ or complex numbers $\C$, and $L^{0}(\mathcal{F}, K)$ (${\bar L}^{0}(\mathcal{F}, \R)$) the algebra of all equivalence classes of ${\mathcal F}$-measurable $K$-valued (accordingly, extended real-valued) random variables on $\Omega$. We write $L^0$ and ${\bar L}^0$ for $L^{0}(\mathcal{F}, \R)$ and ${\bar L}^{0}(\mathcal{F}, \R)$, respectively.

As usual, ${\bar L}^{0}$ is partially ordered by $\xi\leqslant\eta$ iff $\xi^{0}(\omega)\leq\eta^{0}(\omega)$ for $P$-almost all $\omega\in \Omega$, where $\xi^0$ and $\eta^0$ are arbitrarily chosen representatives of $\xi$ and $\eta$, respectively. According to \cite{Dunford}, $({\bar L}^0,\leqslant)$ is a complete lattice, and $( L^0,\leqslant)$ is a conditionally complete lattice. For a subset $H$ of ${\bar L}^0$, $\vee H$ stands for the supremum of $H$, and if $H$ is upward directed, namely there exists $c\in H$ for any $a, b\in H$ such that $a\leqslant c$ and $b\leqslant c$, then there exists a sequence $\{a_n,n\in\N\}$ in $H$ such that $\{a_n,n\in\N\}$ converges to $\vee A$ in a nondecreasing way.

${\tilde I}_A$ always denotes the equivalence class of $I_A$, where $A\in {\mathcal F}$ and $I_A$ is the characteristic function of $A$. For any $\xi\in L^0$, $|\xi|$ denotes the equivalence class of $|\xi^0|: \Omega\to [0,\infty)$ defined by $|\xi^0|(\omega)=|\xi^0(\omega)|$, where $\xi^0$ is an arbitrarily chosen representative of $\xi$.

Denote $L^{0}_{+}=\{\xi\in L^{0}\,|\,\xi\geqslant 0\}$.

Definition \ref{RNModule} is essentially adopted from \cite{Guo-some} excepte for some changes of notations, this form was first employed by \cite{Guo-JFA}.

\begin{definition}\label{RNModule}
An ordered pair $(E,\|\cdot\|)$ is called a random normed module (briefly, an RN module) over $K$ with base $(\Omega,{\mathcal F},P)$ if $E$ is a left module over $L^{0}(\mathcal{F}, K)$, and the mapping $\|\cdot\|: E\to L^0_+$ satisfies:

\noindent(1) $\|x\|=0$ if and only if $x=\theta$( the null element of $E$);\\
\noindent(2)  $\|\xi x\|=|\xi|\|x\|$ for all $\xi\in L^{0}(\mathcal{F}, K)$ and $x\in E$;\\
\noindent(3) $\|x_1+x_2\|\leqslant \|x_1\|+\|x_2\|$ for all $x_1, x_2\in E$.

\end{definition}

In this paper, given an RN module $(E,\|\cdot\|)$, $E$ is always endowed with the $(\varepsilon,\lambda)$-topology. It suffices to say that the $(\varepsilon,\lambda)$-topology is a metrizable linear topology, a sequence $\{x_n,n\in\N\}$ in $E$ converges in the $(\varepsilon,\lambda)$-topology to $x$ iff the sequence $\{\|x_n-x\|,n\in \N\}$ in $L^0_+$ converges in probability to $0$.

Specially, $(L^{0}(\mathcal{F}, K), |\cdot|)$ is an RN module, and the $(\varepsilon,\lambda)$-topology is exactly the topology of convergence in probability.

Let $(E,\|\cdot\|)$ be an RN module over $K$ with base $(\Omega,{\mathcal F},P)$. $E^\ast$ denotes the
$L^{0}(\mathcal{F}, K)$-module of all continuous module homomorphisms $f$ from $(E,\|\cdot\|)$ to $(L^{0}(\mathcal{F}, K), |\cdot|)$. According to \cite[Theorem 3.1]{Guo-ext}, given a linear mapping $f: E\to L^{0}(\mathcal{F}, K)$, then $f\in E^\ast$ if and only if $f$ is a.s bounded, which means that for some $\xi\in L^0_+$, $|f(x)|\leq \xi\|x\|$ holds for all $x\in E$. Define $\|f\|^\ast=\vee\{|f(x)|:x\in E, \|x\|\leqslant 1\}$ for each $f\in E^\ast$, then $(E^\ast, \|\cdot\|^\ast)$ is an RN module, called the random conjugate space of $(E,\|\cdot\|)$.

\section{The Orlicz space generated from a random normed module}
 \label{}
 \subsection{Some basic facts of Orlicz
space theory}

We shortly review some basic facts of Orlicz
space theory. A function $\Phi: [0,\infty)\to [0, \infty]$ is called a Young function if it is convex, left-continuous, $\lim_{t\to 0}\Phi(t)=\Phi(0)=0$, and $\lim_{t\to \infty}\Phi(t)=\infty$. It is easy to see that
$\Phi$ is increasing, and continuous except
possibly at a single point, where it jumps to $\infty$. The conjugate of $\Phi$ defined by
$$\Psi(s):=\sup\limits_{t\geq 0}\{ts-\Phi(t)\},~\quad s\geq 0,$$
is also a Young function. We can check that the conjugate of $\Psi$ is $\Phi$.

In the sequel, we use $E[\xi]$ to denote $\xi$'s expectation with respect to the probability $P$.

The Orlicz space (over $(\Omega,{\mathcal F},P)$) corresponding to $\Phi$ is given by:
$$L^\Phi=\{\xi\in L^0: E[\Phi(c|\xi|)]<\infty~\mbox{for som~} c>0\},$$
and the Orlicz heart is given by:
$$M^\Phi=\{\xi\in L^0: E[\Phi(c|\xi|)]<\infty~\mbox{for all~} c>0\}.$$
The Luxemburg norm
$$|\xi|_{\Phi L}=\inf\{\lambda>0: E[\Phi(\frac{|\xi|}{\lambda})]\leq 1\},$$
and the Orlicz norm
$$|\xi|_{\Phi O}=\sup\{|E[\xi\eta]|: \eta\in L^\Psi, |\eta|_{\Psi L}\leq 1\}$$
are equivalent norms on $L^\Phi$ under which and with the usual partial order $L^\Phi$ is a Banach lattice.

 If $\Phi(t)=\infty$ for some $t\in (0,\infty)$, then $M^\Phi$ is the trivial space $\{0\}$. Thus in the sequel, $\Phi$ is assumed to be real valued, equivalently, $\Phi$ is continuous. According to \cite[Theorem 2.1.14]{ES}, we always have that ${\mathcal S}\subset L^\infty\subset M^\Phi$, where ${\mathcal S}$ stands for the set of all simple measurable functions. Moreover, ${\mathcal S}$ is dense in $(M^\Phi, |\cdot|_{\Phi L})$, which implies that $M^\Phi$ is the $|\cdot|_{\Phi L}$-closure of $L^\infty$ in $L^\Phi$.
 Furthermore, according to \cite[Theorem 2.2.11]{ES}, the norm dual of $(M^\Phi, |\cdot|_{\Phi L})$ ($(M^\Phi, |\cdot|_{\Phi O})$) is given by $(L^\Psi,|\cdot|_{\Psi O})$(accordingly, $(L^\Psi,|\cdot|_{\Psi L})$), where $\eta\in L^\Psi$ is identified with the bounded linear functional $f_\eta: M^\Phi\to {\mathbb R}$ defined by
$f_\eta(\xi)=E[\xi\eta],\forall~\xi\in M^\Phi$.

In the following, we give some simple examples, where we use $|\cdot|_p$ to denote the usual $L^p$ norm.
\begin{example}

1. If $\Phi(t)=t$, then $\Psi(s)=0$ for $s\leq 1$, and $\infty$ otherwise. We have:\\
$M^\Phi=L^\Phi=L^1$, $|\cdot|_{\Phi L}=|\cdot|_{\Phi O}=|\cdot|_1$,\\
$L^\Psi=L^\infty, M^\Psi=\{0\}, |\cdot|_{\Psi O}=|\cdot|_{\Psi L}=|\cdot|_\infty$.

2. If $\Phi(x)=t^p$ for $p\in (1,\infty)$, then $\Psi(s)=p^{1-q}q^{-1}s^q$. We have:\\
$M^\Phi=L^\Phi=L^p$, $|\cdot|_{\Phi L}=|\cdot|_p$, $|\cdot|_{\Phi O}=p^{\frac{1}{p}}q^{\frac{1}{q}}|\cdot|_p$,\\
$M^\Psi=L^\Psi=L^q$, $|\cdot|_{\Psi O}=|\cdot|_q$, $|\cdot|_{\Psi L}=p^{-\frac{1}{p}}q^{-\frac{1}{q}}|\cdot|_q$.
\end{example}

\subsection{The Orlicz space generated from a random normed module}

Let $(E,\|\cdot\|)$ be an RN module with base $(\Omega,{\mathcal F},P)$ and $\Phi$ a Young function. We introduce the Orlicz space corresponding to  $\Phi$ generated from $E$ as:
$$L^\Phi(E)=\{x\in E: \|x\|\in L^\Phi\},$$
and the Orlicz heart of $E$ as: $$M^\Phi(E)=\{x\in E: \|x\|\in M^\Phi\}.$$
Induced by the norm on $L^\Phi$, the Orlicz norm $\|\cdot\|_{\Phi O}$ and the Luxemburg norm $\|\cdot\|_{\Phi L}$ on $L^\Phi(E)$ are given by
$$\|x\|_{\Phi O}=\big|\|x\|\big|_{\Phi O},$$
and
$$\|x\|_{\Phi L}=\big|\|x\|\big|_{\Phi L},$$
for each $x\in L^\Phi(E)$, respectively.

\begin{example}
1. When $(E,\|\cdot\|)=(L^0({\mathcal F},K),|\cdot|)$, $L^\Phi(E)$ and $M^\Phi(E)$ are exactly $L^\Phi$ and $M^\Phi$, respectively.

2.If $\Phi(t)=t^p$ for $p\in [1,\infty)$, then we have $M^\Phi(E)=L^\Phi(E)=L^p(E)$ and $\|\cdot\|_{\Phi L}=\|\cdot\|_p$, where $L^p(E)=\{x\in E: (E[\|x\|^p])<\infty\}$ and $\|x\|_p=(E[\|x\|^p])^\frac{1}{p}$.
Thus $(M^\Phi(E), \|\cdot\|_{\Phi L})$ is exactly the abstract $L^p$-space generated from $E$.

3. Let $\Psi(y)=0$ for $y\leq 1$, and $\infty$ otherwise. Then $M^\Psi(E)=\{0\}, L^\Psi(E)=L^\infty(E):=\{x\in E: \|x\|\in L^\infty\}$ and $\|x\|_{\Psi O}=\|x\|_\infty$.
\end{example}

\begin{proposition}
Let $(E,\|\cdot\|)$ be a complete RN module and $\Phi$ a continuous Young function. Then both $(L^\Phi(E),\|\cdot\|_{\Phi L})$ and $(M^\Phi(E), \|\cdot\|_{\Phi L})$ are Banach spaces.
\end{proposition}
{\em proof.} We only need to prove the completeness. To show the completeness of $(L^\Phi(E),\|\cdot\|_{\Phi L})$, let $\{x_n,n\in \N\}$ be an arbitrary Cauchy sequence in $(L^\Phi(E), \|\cdot\|_{\Phi L})$. Then $\{x_n,n\in \N\}$ must also be a Cauchy sequence in $(E,\|\cdot\|)$, otherwise, there exist $\varepsilon>0$ and $\lambda\in (0,1)$ such that for any given $N_0\in\N$, we can always find some $m,n\geq N_0$ such that $P\{\omega\in\Omega: \|x_m-x_n\|(\omega)\geq \varepsilon\}\geq \lambda$, then  $$1\geq E[\Phi(\frac{\|x_m-x_n\|}{c})]\geq \lambda\Phi(\frac{\varepsilon}{c})$$
 yields that $c\geq c_0=\varepsilon(\Phi^{-1}({1\over\lambda}))^{-1}$, which means that $\|x_m-x_n\|_{\Phi L}\geq c_0>0$ for these m's and n's, contradicting to the assumption that $\{x_n,n\in \N\}$ is a Cauchy sequence in $(L^\Phi(E), \|\cdot\|_{\Phi L})$. Then using the fact that $(E,\|\cdot\|)$ is complete, $\{x_n,n\in \N\}$ converges to some $x\in E$, namely, $\|x_n-x\|\to 0$ in probability as $n\to \infty$. Using Fatou's lemma, we get $$\lim_{n\to\infty}E[\Phi(\frac{\|x_n-x\|}{c})]\leq \lim_{n\to\infty}\lim_{m\to\infty}E[\Phi(\frac{\|x_n-x_m\|}{c})]\leq 1, \forall c>0,$$
which implies that $\lim\limits_{n\to\infty}\|x_n-x\|_{\Phi L}=0$.

$(M^\Phi(E), \|\cdot\|_{\Phi L})$ is a Banach space follows the next proposition.
\hfill\done

\begin{proposition}\label{dense}
Let $(E,\|\cdot\|)$ be an RN module and $\Phi$ a continuous Young function, then $M^\Phi(E)$ is the $\|\cdot\|_{\Phi L}$-closure of $L^\infty(E)$ in $L^\Phi(E)$.
\end{proposition}
{\em proof.} Let $\{x_n,n\in \N\}$ be a sequence in $L^\infty(E)$ and $x\in L^\Phi(E)$ such that $\|x_n-x\|_{\Phi L}\to 0$ as $n\to\infty$. Since $|\|x_n\|-\|x\||\leq \|x_n-x\|$, we have $|(\|x_n\|-\|x\|)|_{\Phi L}\to 0$ as $n\to\infty$. Then $\|x\|\in M^\Phi$ follows $\|x_n\|\in L^\infty,\forall n\in\N$, thus $x\in M^\Phi(E)$, which means that the $\|\cdot\|_{\Phi L}$-closure of $L^\infty(E)$ is contained in $M^\Phi(E)$.

Fixed $x\in M^\Phi(E)$, for each $n\in\N$, let $A_n=\{\omega\in\Omega:\|x\|^0(\omega)\leq n\}$ and take $x_n={\tilde I}_{A_n}x$, where $\|x\|^0$ is an arbitrarily chosen representative of $\|x\|$, then $x_n\in L^\infty(E)$ and $\|x_n-x\|_{\Phi L}=|{\tilde I}_{A^c_n}\|x\||_{\Phi L}\to 0$ as $n\to\infty$. Thus the $\|\cdot\|_{\Phi L}$-closure of $L^\infty(E)$ contains $M^\Phi(E)$.
\hfill\done

\begin{remark}
Let $(E,\|\cdot\|)$ be a complete RN module and $\Phi$ a continuous Young function. Since the Orlicz norm $\|\cdot\|_{\Phi O}$ and the Luxemburg norm $\|\cdot\|_{\Phi L}$ are equivalent norms on $L^\Phi(E)$, both $(L^\Phi(E), \|\cdot\|_{\Phi O})$ and $(M^\Phi(E), \|\cdot\|_{\Phi O})$ are also Banach spaces.
\end{remark}

\section{Dual representation of the conjugate space of $M^{\Phi}(E)$}

We first state the main result as follows.

\begin{theorem}\label{main}
Let $(E,\|\cdot\|)$ be an RN module over $K$ with base $(\Omega,{\mathcal F},P)$, $(E^\ast,\|\cdot\|^\ast)$ its random conjugate space, and $\Phi$ a continuous Young functions with conjugate $\Psi$. Then $$(M^\Phi(E),\|\cdot\|_{\Phi L})^\prime\cong (L^\Psi(E^\ast),\|\cdot\|_{\Psi O}),$$ where the isometric isomorphism $T:(L^\Psi(E^\ast),\|\cdot\|_{\Psi O})\to (M^\Phi(E),\|\cdot\|_{\Phi L})^\prime$ is given by $$[Tf](x)=E[f(x)], \forall x\in M^\Phi(E),$$
for each $f\in L^\Psi(E^\ast)$.
\end{theorem}

For the sake of clearness, the proof of Theorem \ref{main} is divided into the following two Lemmas
\ref{isometric} and \ref{onto}. Lemma \ref{isometric} shows that $T$ is well defined, and isometric and Lemma \ref{onto} shows that $T$ is surjective.

In the following, $U(E):=\{x\in E: \|x\|\leqslant 1\}$ denotes the random closed unit ball of $E$.

\begin{lemma}\label{isometric}
$T$ is well defined and isometric.
\end{lemma}
{\em proof.}
  For any fixed $f\in L^\Psi(E^\ast)$, we will prove $Tf\in (M^\Phi(E),\|\cdot\|_{\Phi L})^\prime$ and $\|Tf\|=\|f\|_{\Psi O}$.

For any $x\in M^\Phi(E)$, according to ``H\"{o}lder inequality'', we have: $$|[Tf](x)|=|E[f(x)]|\leq E[\|f\|^\ast\|x\|]\leq |\|f\|^\ast|_{\Psi O} |\|x\||_{\Phi L}=\|f\|_{\Psi O} \|x\|_{\Phi L}.$$
This shows that $Tf\in (M^\Phi(E),\|\cdot\|_{\Phi L})^\prime$ and $\|Tf\|\leq \|f\|_{\Psi O}$, namely, $T$ is well-defined. We remain to show that $\|Tf\|=\|f\|_{\Psi O}$, or equivalently, $\|Tf\|\geq \|f\|_{\Psi O}$.

Note that\begin{eqnarray*}
           \|f\|_{\Psi O}=|\|f\|^\ast|_{\Psi O}&=& \sup\{|E[\|f\|^\ast\xi]|:\xi\in M^\Phi, |\xi|_{\Phi L}\leq 1\}\\
           &=& \sup\{E[\|f\|^\ast\xi]:\xi\in M^\Phi, \xi\geq 0, |\xi|_{\Phi L}\leq 1\},
         \end{eqnarray*}
thus for any fixed $\xi\in M^\Phi$ with $\xi\geq 0, |\xi|_{\Phi L}\leq 1$, it suffices to show $\|Tf\|\geq E[\|f\|^\ast\xi]$. It is easy to verify that the family $\{|f(x)|: x\in U(E)\}$ is upward directed, thus there exists a sequence $\{x_n,n\in\mathbb{N}\}$ in $U(E)$ such that $\{|f(x_n)|, n\in \mathbb{N}\}$ converges to $\vee\{|f(x)|: x\in U(E)\}=\|f\|^\ast$ in a nondecreasing way. Further, we can assume that $f(x_n)=|f(x_n)|$ for every $n$, otherwise we can replace each $x_n$ with $(sgn f(x_n))x_n$ ( here $sgn(z)$ for an element $z\in L^0({\mathcal F}, K)$ means the equivalence class of $sgn(z^0)$ defined by $sgn(z^0)(\omega)=\frac{|z^0(\omega)|}{z^0(\omega)}$ if $z^0(\omega)\neq 0$, and $0$ otherwise, where $z^0$ is an arbitrarily chosen representative of $z$ ). Thus, $\lim_{n\to\infty}f(\xi x_n)=\lim_{n\to\infty}\xi f(x_n)=\xi\|f\|^\ast$.  For each $n$, since $\|\xi x_n\|=\xi\|x_n\|\leqslant \xi$, we have $\|\xi x_n\|_{\Phi L}\leq |\xi|_{\Phi L}\leq 1$, thus $E[f(\xi x_n)]=[Tf](\xi x_n)\leq \|Tf\|$. Then according to Levi's monotone convergence theorem, we finally get $$E[\|f\|^\ast\xi]=\lim_{n\to\infty}E[f(\xi x_n)]\leq \|Tf\|,$$
which completes the proof.
\hfill\done

\begin{lemma}\label{onto}
$T$ is surjective.
\end{lemma}
{\em proof}
Let $F$ be an arbitrary element in $(M^\Phi(E),\|\cdot\|_{\Phi L})^\prime$, we want to prove that there exists an $f\in L^\Psi(E^\ast)$ such that $F=Tf$.

For any fixed $x\in M^\Phi(E)$, define $\mu_x: {\mathcal F}\to K$ by $\mu_x(A)=F({\tilde I}_Ax), \forall A\in {\mathcal F}$, then $\mu_x$ is a countably additive $K$-valued measure, which is absolutely continuous with respect to the probability measure $P$. Thus according to Radon-Nikod\'{y}m theorem, there exists an unique $\xi_x\in L^1$ such that $\mu_x(A)=F({\tilde I}_Ax)=E[{\tilde I}_A\xi_x], \forall A\in {\mathcal F}$. Moreover, $|\mu_x|(\Omega)=E[|\xi_x|]\leq \|F\|\|x\|_{\Phi L}$, since $|\mu_x(A)|=|F({\tilde I}_Ax)|\leq \|F\|\|x\|_{\Phi L}, \forall A\in {\mathcal F}$. Define $g: (M^\Phi(E),\|\cdot\|_{\Phi L}) \to (L^1,|\cdot|_1)$ by $g(x)=\xi_x, \forall x\in M^\Phi(E)$, then $g$ is a bounded linear operator, and $g({\tilde I}_Ax)={\tilde I}_Ag(x),\forall A\in {\mathcal F}, x\in M^\Phi(E)$. Immediately, we have that for each $x\in M^\Phi(E)$, $g(\xi x)=\xi g(x)$ holds for every simple function $\xi\in L^{0}(\mathcal{F}, K)$. Further we verify that $g(\xi x)=\xi g(x)$ holds for every $x\in U(E)$ and $\xi\in M^\Phi$. In fact, fix $x\in U(E)$ and $\xi\in M^\Phi$, according to \cite[Theorem 2.1.14]{ES}, there exists a sequence $\{\xi_n,n\in \mathbb{N}\}$ consisting of simple measurable functions such that $|\xi_n-\xi|_{\Phi L}\to 0$ which implies that $\xi_n\to \xi$ in probability and $\|\xi_n x-\xi x\|_{\Phi L}=|(\xi_n-\xi)\|x\||_{\Phi L}\leq |\xi_n-\xi|_{\Phi L}\to 0$ as $n\to\infty$, thus $g(\xi x)=L^1-\lim_{n\to\infty}g(\xi_n x)=L^1-\lim_{n\to\infty} \xi_n g(x)=\xi g(x)$.

Consider the subset $\{|g(x)|: x\in U(E)\}$ of $L^1\subset L^0$, since $|g({\tilde I}_A x+{\tilde I}_{A^c}y)|={\tilde I}_A|g(x)|+{\tilde I}_{A^c}|g(y)|$ holds for every $x,y\in U(E)$ and $A\in {\mathcal F}$, we see that $\{|g(x)|: x\in U(E)\}$ is upward directed. As in the proof of Lemma \ref{isometric}, there exists a sequence $\{x_n,n\in\mathbb{N}\}$ in $U(E)$ such that $|g(x_n)|=g(x_n)$ converges to $X_g:=\vee\{|g(x)|: x\in U(E)\}\in {\bar L}^0_+$ in a nondecreasing way as $n\to\infty$. Then, for any $\xi\in M^\Phi$ with $\xi\geq 0$, by Levi's monotone convergence theorem, $$E[\xi X_g]=\lim_{n\to\infty}E[\xi g(x_n)]=\lim_{n\to\infty}E[g(\xi x_n)]\leq \lim_{n\to\infty}\|F\|\|\xi x_n\|_{\Phi L}\leq \|F\||\xi|_{\Phi L}<\infty,$$
which implies that $X_g\in L^\Psi$ with $|X_g|_{\Psi O}\leq \|F\|$.

Define $f: E\to L^{0}(\mathcal{F}, K)$ by $f(x)=\lim_{n\to \infty}g({\tilde I}_{A_n}x)$ for each $x\in E$, where $A_n$ is taken as in the proof of Proposition \ref{dense} and the limit on the right side is taken with respect to the convergent in probability. Since $|g({\tilde I}_{A_m}x)-g({\tilde I}_{A_n}x)|=|g(({\tilde I}_{A_m}-{\tilde I}_{A_n})x)|\leqslant X_g|{\tilde I}_{A_m}-{\tilde I}_{A_n}|\|x\|\to 0$ in probability as $m,n\to\infty$, we see that $f$ is well defined. Moreover, it is easily seen that $f$ is linear, $f(x)=g(x)$ holds for every $x\in L^\infty (E)$ and $|f(x)|\leqslant X_g\|x\|, \forall x\in E$, it follows from \cite[Lemma 2.12]{Guo-JFA} that $f$ is $L^{0}(\mathcal{F}, K)$-linear, thus $f\in E^\ast$ and $$\|f\|^\ast=\vee\{|f(x)|:x\in U(E)\}=\vee\{|g(x)|:x\in U(E)\}=X_g,$$
which implies that $f\in L^\Psi(E^\ast)$.

We remain to show $F=Tf$. Note that $$[Tf](x)=E[f(x)]=E[g(x)]=F(x), \forall x\in L^\infty(E),$$
since $L^\infty(E)$ is a dense subset of $(M^\Phi(E),\|\cdot\|_{\Phi L})$ by Proposition \ref{dense}, the two continuous functionals $F$ and $Tf$ must equal to each other on the whole space $(M^\Phi(E),\|\cdot\|_{\Phi L})$.
\hfill\done

Obviously, if the Luxemburg norm $\|\cdot\|_{\Phi L}$ on $M^\Phi(E)$ is replaced by Orlicz norm $\|\cdot\|_{\Phi O}$, then the operator norm on the dual space $L^\Psi(E)$ changes accordingly from Orlicz norm $\|\cdot\|_{\Psi O}$ to Luxemburg norm $\|\cdot\|_{\Psi L}$. Precisely, we have:

\begin{proposition}\label{varmain}
 Let $(E,\|\cdot\|)$ be an RN module over $K$ with base $(\Omega,{\mathcal F},P)$, $(E^\ast,\|\cdot\|^\ast)$ its random conjugate space, and $\Phi$ a continuous Young functions with conjugate $\Psi$. Then $$(M^\Phi(E),\|\cdot\|_{\Phi O})^\prime\cong (L^\Psi(E^\ast),\|\cdot\|_{\Psi L}),$$ where the isometric isomorphism is the same as Theorem \ref{main}.
\end{proposition}

Now we add some conditions on $\Phi$ and $\Psi$. A Young function $\Phi: [0,+\infty)\to [0,+\infty)$ is said to satisfy the $\triangle_2$- condition, denoted by $\Phi\in \triangle_2$, if there is a constant $k>0$ such that $\Phi(2u)\leq k\Phi(u)$ for every $u\in [0,+\infty)$. When a Young function $\Phi$ satisfies the $\triangle_2$- condition, the Olicz space $L^\Phi$ and the Orlicz heart $M^\Phi$ coincide (see \cite[Theorem 2.1.17]{ES}). Assume that $\Phi$ and $\Psi$ is a pair of conjugate Young functions such that $\Phi\in \triangle_2$ and $\Psi\in \triangle_2$, then we have $L^\Phi(E)=M^\Phi(E)$ and $L^\Psi(E^\ast)=M^\Psi(E^\ast)$. Using Theorem \ref{main} and Proposition \ref{varmain}, we obtain the following relations: $$(L^\Phi(E),\|\cdot\|_{\Phi L})^{\prime\prime}\cong (L^\Psi(E^\ast),\|\cdot\|_{\Psi O})^\prime\cong (L^\Phi(E^{\ast\ast}),\|\cdot\|_{\Phi L}).$$

Recall that an RN module $E$ is said to be random reflexive if the canonical embedding mapping $J: E\to E^{\ast\ast}$ defined by $[J(x)](f)=f(x), \forall x\in E. f\in E^\ast$, is surjective. It is straightforward to show the following:

\begin{proposition}\label{reflex}
 Let $(E,\|\cdot\|)$ be an RN module over $K$ with base $(\Omega,{\mathcal F},P)$, $\Phi$ and $\Psi$ a pair of conjugate Young functions such that $\Phi\in \triangle_2$ and $\Psi\in \triangle_2$. Then the following statements are equivalent:\\
  (1) $(E,\|\cdot\|)$ is random reflexive;\\
  (2) $(L^\Phi(E),\|\cdot\|_{\Phi L})$ is reflexive;\\
  (3) $(L^\Phi(E),\|\cdot\|_{\Phi O})$ is reflexive.
\end{proposition}

\begin{remark}
In Theorem \ref{main}, if we choose $\Phi(t)=t^p$ for $p\in [1,\infty)$, then we get $$(L^p(E),\|\cdot\|_p)^\prime\cong (L^q(E^\ast),\|\cdot\|_q). $$
In Proposition \ref{reflex}, if we choose $\Phi(t)=t^p$ for $p\in (1,\infty)$, then we have that $(E,\|\cdot\|)$ is random reflexive if and only if $(L^p(E),\|\cdot\|_p)$ is reflexive. Thus the results in this section generalize some known results established by Guo \cite{Guo-Xiamen,Guo-rep}.
\end{remark}

\section{Strict convexity and uniform convexity of $L^{\Phi}(E)$}
\subsection{Main results}
It is well known that strictly convex and uniformly convex Banach spaces have played key roles in many important topics in nonlinear functional analysis and geometry of Banach spaces. This section is devoted to study the strict convexity and uniform convexity of $L^{\Phi}(E)$. We will establish a basic connection between uniform convexity (strict convexity) of $L^{\Phi}(E)$ and random uniform convexity (accordingly, random strict convexity) of $E$ together with uniform convexity (accordingly, strict convexity) of $L^{\Phi}$.

Considering that strict convexity and uniform convexity are properties of norm, and as introduced as above, there are two norms, namely Luxemburg norm $\|\cdot\|_{\Phi L}$ and the Orlicz norm $\|\cdot\|_{\Phi O}$ on the space $L^{\Phi}(E)$, we therefore need to study the strict convexity and uniform convexity of $L^{\Phi}(E)$ under the two norms, respectively. Surprisingly, it turns out that we can discuss the two cases in a unified way rather than separately. In the sequel, the norm on $L^\Phi$ is denoted by $|\cdot|_\Phi$, which can be either the Luxemburg norm $|\cdot|_{\Phi L}$ or the Orlicz norm $|\cdot|_{\Phi O}$, and accordingly, the norm on $L^\Phi(E)$ induced by $|\cdot|_\Phi$ is denoted by $\|\cdot\|_\Phi$.

The main results in this section are as follows:
\begin{theorem}\label{strcon}
Let $(E,\|\cdot\|)$ be a complete RN module with base $(\Omega,{\mathcal F},P)$ such that $E$ has full support, and $\Phi$ a given Young function, then $(L^\Phi(E),\|\cdot\|_\Phi)$ is strictly convex if and only if $(L^\Phi, |\cdot|_\Phi)$ is strictly convex and $(E,\|\cdot\|)$ is random strictly convex.
\end{theorem}

\begin{theorem}\label{unicon}
Let $(E,\|\cdot\|)$ be a complete RN module with base $(\Omega,{\mathcal F},P)$ such that $E$ has full support, and $\Phi$ a given Young function, then $(L^\Phi(E),\|\cdot\|_\Phi)$ is uniformly convex if and only if $(L^\Phi, |\cdot|_\Phi)$ is uniformly convex and $(E,\|\cdot\|)$ is random uniformly convex.
\end{theorem}

In the remainder of this section, we first recall some notions and notations involved in the two theorems above in Subsection \ref{NAN}, and then give the proof of the Theorem \ref{strcon} in Subsection \ref{STR} and the proof of Theorem \ref{unicon} in Subsection \ref{UNI}.

\subsection{Some notions and notations}\label{NAN}

Recall that a normed space $(X,\|\cdot\|)$ is said to be:\\
Strictly convex: if for any two distinct elements $x,y\in X$ with $\|x\|=\|y\|=1$, we have $\|\frac{x+y}{2}\|<1$;\\
Uniformly convex: if for every $\epsilon\in (0,2]$, there exists $\delta(\epsilon)>0$ such that for any $x,y\in X$ such that $\|x\|=\|y\|=1$ and $\|x-y\|\geq \epsilon$, we have $\|\frac{x+y}{2}\|\leq 1-\delta(\epsilon)$.

We then recall the notions of random strict convexity and random uniform convexity of an RN module which were introduced by Guo and Zeng \cite{GZ1}.

Let $A\in {\mathcal F}$, then the equivalence class of $A$, denoted by $\tilde A$, is defined by $\tilde A=\{B\in {\mathcal F}: P(A\triangle B)=0\}$, where $A\triangle B=(A\setminus B)\cup (B\setminus A)$ is the symmetric difference of $A$ and $B$, and $P(\tilde A)$ and $I_{\tilde A}$ are defined to be $P(A)$ and ${\tilde I_A}$, respectively. Denote ${\tilde {\mathcal F}}=\{\tilde A: A\in {\mathcal F}\}$. For two ${\mathcal F}$-measurable sets $G$ and $D$, $G\subset D$ a.s. means $P(G\setminus D)=0$, in which case we also say $\tilde G\subset \tilde D$; $\tilde G\cap\tilde D$ denotes the equivalence class determined by $G\cap D$. For any $\xi,\eta\in L^0$ and $A\in {\mathcal F}$, $\xi>\eta $ on $\tilde A$ means that $\xi^0(\omega)>\eta^0(\omega)$ for almost all $\omega\in A$, and $[~\xi\geq\eta~]$ means the equivalence class of ${\mathcal F}$-measurable set $\{\omega\in \Omega: \xi^0(\omega)\geq \eta^0(\omega)\}$, where $\xi^0$ and $\eta^0$ are arbitrarily chosen representatives of $\xi$ and $\eta$, respectively. Other similar symbols are easily understood in an analogous manner.

Let $(E,\|\cdot\|)$ be an RN module. For any $x, y\in E$, denote $[~\|x\|\neq 0~]$ by $A_x$, called the support of $x$, and let $B_{xy}=A_x\cap A_y\cap A_{x-y}$. The random unite sphere of $E$ refers to $$S(E)=\{x\in E: x\neq \theta, ~\|I_{A_x}x\|=I_{A_x}\}.$$

Further we employ the following notations:
\begin{align*}
& L_\#=\{\xi\in L^0_+: \exists \lambda\in {\mathbb R}, \lambda>0 \text{~such that~} \xi\geqslant \lambda \};\\
& L_\#[0,1]=\{\xi\in L_\#: \xi\leqslant 1\}; \\
& L_\#[0,2]=\{\xi\in L_\#: \xi\leqslant 2\}.
\end{align*}

\begin{definition}(\cite{GZ1}).
 Let $(E,\|\cdot\|)$ be an RN module with base $(\Omega,{\mathcal F},P)$. $E$ is said to be random strictly convex, if for any $x,y\in S(E)$ with $P(B_{xy})>0$, we have $\|\frac{x+y}{2}\|<1$ on $B_{xy}$; and $E$ is said to be random uniformly convex if for any $\epsilon\in L_\#[0,2]$ there exists a $\delta\in L_\#[0,1]$ such that the following condition holds: $I_D\|x-y\|\geqslant \epsilon I_D$ implies that $I_D\|\frac{x+y}{2}\|\leqslant (1-\delta)I_D$ for any $x, y\in S(E)$ and $D\in {\tilde{\mathcal F}}$ with $D\subset B_{xy}$ and $P(D)>0$.
\end{definition}

Let $(E,\|\cdot\|)$ be an RN module, set $\xi(E)=\vee\{\|x\|: x\in E\}$, and denote $H(E)=[~\xi(E)\neq 0]\in {\tilde {\mathcal F}}$, namely, $H(E)$ is the esssup of the family $\{[~\|x\|\neq 0~]: x\in E\}$, called the support of $E$. When $H(E)=\tilde \Omega$, $E$ is called having full support. We restate \cite[Lemma 3.1]{BBKS} as follows.

\begin{proposition}\label{fulspp}\
Assume that $(E,\|\cdot\|)$ is a complete RN module, then there exists an $x_0\in E$ such that $\|x_0\|=I_{H(E)}$, specially, when $(E,\|\cdot\|)$ has full support, then there exists an $x_0\in E$ such that $\|x_0\|=1$.
\end{proposition}

\begin{remark}
  Let $E$ be an RN module which is nontrivial, namely $P(H(E))>0$. If $E$ does not have full support, then on ${\tilde \Omega}\setminus H(E)$, we have $x=0,\forall x\in E$, which means that the part ${\tilde \Omega}\setminus H(E)$ is redundant for $E$. In this case, to capture the essential properties of $E$, we may regard $E$ as an RN module over the smaller probability space $(H, {\mathcal F}\cap H, P_H)$, where $H\in {\mathcal F}$ is an arbitrarily representative of $H(E)\in \tilde{\mathcal F}$ and $P_H(A\cap H)=\frac{P(A\cap H)}{P(H)}, \forall A\in {\mathcal F}$, then $E$ has full support. Thus the assumption on $E$ in Theorem \ref{strcon} and Theorem \ref{unicon} is not a substantial restriction.
\end{remark}

\subsection{Strict convexity of $L^{\Phi}(E)$}\label{STR}

In this subsection, we prove Theorem \ref{strcon}. To clarify the implications of the assumption on $E$ in Theorem \ref{strcon}, we divide Theorem \ref{strcon} into two propositions--Proposition \ref{strsuf} and Proposition \ref{strnec}.

\begin{lemma}\label{strlem}
If $(L^\Phi, |\cdot|_\Phi)$ is strictly convex, then for any two distinct elements $\xi,\eta\in L^\Phi$ with $\xi\geqslant \eta\geqslant 0$, it holds that $|\xi|_\Phi>|\eta|_\Phi$.
\end{lemma}

{\em proof.}
The assumption $\xi\geqslant \eta\geqslant 0$ yields that $|\xi|_\Phi\geq |\eta|_\Phi$. We prove the conclusion by contradiction. Suppose that $|\xi|_\Phi=|\eta|_\Phi=\lambda$. Since $\xi\neq \eta$, $\lambda$ must be a positive number. Let $\xi_0=\frac{\xi}{\lambda}$ and $\eta_0=\frac{\eta}{\lambda}$, then $|\xi_0|_\Phi=|\eta_0|_\Phi=1$, clearly $\xi_0\geqslant\frac{\xi_0+\eta_0}{2}\geqslant \eta_0\geqslant 0$, implying that $|\frac{\xi_0+\eta_0}{2}|_\Phi=1$, but this contradicts to the assumption that $(L^\Phi, |\cdot|_\Phi)$ is strictly convex.  \hfill\done

\begin{proposition}\label{strsuf}
 Let $(E,\|\cdot\|)$ be an RN module with base $(\Omega,{\mathcal F},P)$ and $\Phi$ a given Young function, if $(L^\Phi, |\cdot|_\Phi)$ is strictly convex and $(E,\|\cdot\|)$ is random strictly convex, then $(L^\Phi(E),\|\cdot\|_\Phi)$ is strictly convex.
\end{proposition}
{\em proof.}
For any two distinct elements $x,y\in L^\Phi(E)$ with $\|x\|_\Phi=\|y\|_\Phi=1$, it is divided into two cases to show that $\|\frac{x+y}{2}\|_\Phi<1$.

Case 1: When $\|x\|\neq \|y\|$, due to the assumption that $(L^\Phi,|\cdot|_\Phi)$ is strictly convex and $\|x\|_\Phi=\|y\|_\Phi=1$, we have $|\frac{\|x\|+\|y\|}{2}|_\Phi<1$, then $\|\frac{x+y}{2}\|_\Phi<1$ follows immediately from that $\|\frac{x+y}{2}\|\leqslant \frac{\|x\|+\|y\|}{2}$.

Case 2: When $\|x\|=\|y\|$, then $\|\frac{x+y}{2}\|\leqslant \frac{\|x\|+\|y\|}{2}=\|x\|$, and due to the assumption that $(E,\|\cdot\|)$ is random strictly convex  we have $\|\frac{x+y}{2}\|\neq \|x\|$. According to Lemma \ref{strlem}, $\|\frac{x+y}{2}\|_\Phi <\|x\|_\Phi=1$.
\hfill\done

\begin{proposition}\label{strnec}
 Let $(E,\|\cdot\|)$ be a complete RN module with base $(\Omega,{\mathcal F},P)$ such that $E$ has full support and $\Phi$ a given Young function, if $(L^\Phi(E),\|\cdot\|_\Phi)$ is strictly convex, then $(L^\Phi, |\cdot|_\Phi)$ is strictly convex and $(E,\|\cdot\|)$ is random strictly convex.
\end{proposition}

{\em proof.}
For any two distinct elements $\xi,\eta\in L^\Phi$ with $|\xi|_\Phi=|\eta|_\Phi=1$, according to Proposition \ref{fulspp}, there exists an $x_0\in E$ such that $\|x_0\|=1$. Choose $x=\xi x_0$ and $y=\eta x_0$, then $x,y\in L^\Phi(E)$ and $\|x\|_\Phi=\|y\|_\Phi=1$. Note that $x\neq y$, according to the assumption that $(L^\Phi(E),\|\cdot\|_\Phi)$ is strictly convex, we have $|\frac{\xi+\eta}{2}|_\Phi=\|\frac{x+y}{2}\|_\Phi<1$, which means that $(L^\Phi, |\cdot|_\Phi)$ is strictly convex.

We show the random strict convexity of $(E,\|\cdot\|)$ by contradiction. Suppose that there exist $x,y\in S(E)$ with $P(B_{xy})>0$, such that for some $D\in {\tilde{\mathcal F}}$, $D\subset B_{xy}$ with $P(D)>0$, it holds that $I_D\|\frac{x+y}{2}\|=I_D$. Clearly, $\|I_D x\|=I_D\|x\|=I_D=\|I_D y\|$. Set $\lambda=|I_D|_\Phi>0$, then $\|I_D x\|_\Phi=\|I_Dy\|_\Phi=\lambda$. Let $x^\prime=\frac{I_D x}{\lambda}$ and $y^\prime=\frac{I_D y}{\lambda}$, then we have $x^\prime,y^\prime\in L^\Phi(E)$ and $\|x^\prime\|_\Phi=\|y^\prime\|_\Phi=1$. It follows from $D\subset B_{xy}$ that $x^\prime\neq y^\prime$, thus the assumption $(L^\Phi(E),\|\cdot\|_\Phi)$ is strictly convex yields that $\|\frac{x^\prime+y^\prime}{2}\|_\Phi< 1$, equivalently $\|I_D\frac{x+y}{2}\|_\Phi<\lambda$, however $I_D\|\frac{x+y}{2}\|=I_D$ implies that $\|I_D\frac{x+y}{2}\|_\Phi=|I_D|_\Phi=\lambda$.
\hfill\done

Notice that the assumption ``$(E,\|\cdot\|)$ is complete and has full support '' is used only in the process to show that ``$(L^\Phi(E),\|\cdot\|_\Phi)$ is strictly convex'' implies that ``$(L^\Phi, |\cdot|_\Phi)$ is strictly convex'', thus by removing this assumption and adding the assumption of strict convexity of $(L^\Phi, |\cdot|_\Phi)$ we can obtain the following:

\begin{corollary}
Let $(E,\|\cdot\|)$ be an RN module with base $(\Omega,{\mathcal F},P)$ and $\Phi$ a given Young function. Assume that $(L^\Phi, |\cdot|_\Phi)$ is strictly convex, then $(E,\|\cdot\|)$ is random strictly convex if and only if $(L^\Phi(E),\|\cdot\|_\Phi)$ is strictly convex.
\end{corollary}

Since for every $p\in (1,\infty)$, $(L^p, |\cdot|_p)$ is strictly convex, thus we have the following:

\begin{corollary}
Let $(E,\|\cdot\|)$ be an RN module and $1<p<\infty$, then $(E,\|\cdot\|)$ is random strictly convex if and only if $(L^p(E), \|\cdot\|_p)$ is strictly convex.
\end{corollary}

This corollary is exactly Theorem 3.3 in  Guo and Zeng \cite{GZ1}.

\subsection{Uniform convexity of $L^{\Phi}(E)$}\label{UNI}
In this subsection, we prove Theorem \ref{unicon}. As in Subsection \ref{STR}, we divide Theorem \ref{unicon} into two propositions--Proposition \ref{uninec} and Proposition \ref{unisuf}.

\begin{proposition}\label{uninec}
 Let $(E,\|\cdot\|)$ be a complete RN module with base $(\Omega,{\mathcal F},P)$ such that $E$ has full support and $\Phi$ a given Young function. If $(L^\Phi(E),\|\cdot\|_\Phi)$ is uniformly convex, then $(L^\Phi, |\cdot|_\Phi)$ is uniformly convex and $(E,\|\cdot\|)$ is random uniformly convex.
\end{proposition}

{\em proof.}
(1)We first show that $(L^\Phi, |\cdot|_\Phi)$ is uniformly convex. Given $\epsilon\in (0,2]$. For any two elements $\xi,\eta\in L^\Phi$ with $|\xi|_\Phi=|\eta|_\Phi=1$ and $|\xi-\eta|_\Phi\geq \epsilon$, according to Proposition \ref{fulspp}, there exists an $x_0\in E$ such that $\|x_0\|=1$, and if we take $x=\xi x_0$ and $y=\eta x_0$, then we have $\|x\|_\Phi=|\xi|_\Phi=1$, $\|y\|_\Phi=|\eta|_\Phi=1$ and $\|x-y\|_\Phi=|\xi-\eta|_\Phi\geq \epsilon$. Since $(L^\Phi(E),\|\cdot\|_\Phi)$ is uniformly convex, there exists a $\delta\in (0,1]$ such that $|\frac{\xi+\eta}{2}|_\Phi=\|\frac{x+y}{2}\|_\Phi\leq (1-\delta)$. This $\delta$ is decided by $(L^\Phi(E),\|\cdot\|_\Phi)$, not depending on $x, y$ and therefore neither on $\xi,\eta$. This means that $(L^\Phi, |\cdot|_\Phi)$ is uniformly convex.

(2) We show that $(E,\|\cdot\|)$ is random uniformly convex by contradiction. Suppose that $E$ is not random uniformly convex. Then we can find an $\epsilon \in L_\#[0,2]$ such that for any $\delta\in L_\#[0,1]$, there exist $x_\delta, y_\delta\in S(E)$ and $D_\delta\in {\tilde {\mathcal F}}$ satisfying $D_\delta\subset B_{x_\delta y_\delta}$, $P(D_\delta)>0$ and $I_{D_\delta}\|x_\delta-y_\delta\|\geqslant \epsilon I_{D_\delta}$ and $I_{D_\delta}\|\frac{x_\delta+y_\delta}{2}\|>(1-\delta)I_{D_\delta}$ on $D_\delta$. Since $\epsilon \in L_\#[0,2]$, there exists some positive number $\lambda$ such that $\lambda<\epsilon\leq 2$ on $\tilde \Omega$. For this $\lambda$, by the uniform convexity of $L^\Phi(E)$, there exists a number $\delta_1\in (0,1]$ such that for any $u,v\in L^\Phi(E)$ with $\|u\|_\Phi=\|v\|_\Phi=1$, $\|u+v\|_\Phi>2(1-\delta_1)$ implies $\|u-v\|_\Phi<\lambda$. Since $\delta_1$ can be regarded as an element in $L_\#[0,1]$, there exist $x_{\delta_1}, y_{\delta_1}\in S(E)$ and $D_{\delta_1}\in {\tilde {\mathcal F}}$ such that $D_{\delta_1}\subset B_{x_{\delta_1} y_{\delta_1}}$,  $P(D_{\delta_1})>0$ and $I_{D_{\delta_1}}\|x_{\delta_1}-y_{\delta_1}\|\geqslant \epsilon I_{D_{\delta_1}}$, $I_{D_{\delta_1}}\|\frac{x_{\delta_1}+y_{\delta_1}}{2}\|>(1-{\delta_1})I_{D_{\delta_1}}$ on $D_{\delta_1}$.
Let $c=|I_{D_{\delta_1}}|_\Phi >0$, since $\|I_{D_{\delta_1}}x_{\delta_1}\|=\|I_{D_{\delta_1}}y_{\delta_1}\|=I_{D_{\delta_1}}$, we have $\|I_{D_{\delta_1}}x_{\delta_1}\|_\Phi=\|I_{D_{\delta_1}}y_{\delta_1}\|_\Phi=c$. Choose $x=\frac{I_{D_{\delta_1}}x_{\delta_1}}{c}$ and $y=\frac{I_{D_{\delta_1}}y_{\delta_1}}{c}$, then $\|x\|_\Phi=\|y\|_\Phi=1$. From $I_{D_{\delta_1}}\|x_{\delta_1}-y_{\delta_1}\|\geqslant \epsilon I_{D_{\delta_1}}\geqslant \lambda I_{D_{\delta_1}}$ we obtain $c\|x-y\|_\Phi\geq \lambda |I_{D_{\delta_1}}|_\Phi=\lambda c$, and from $I_{D_{\delta_1}}\|\frac{x_{\delta_1}+y_{\delta_1}}{2}\|>(1-{\delta_1})I_{D_{\delta_1}}$ on $D_{\delta_1}$ and Lemma \ref{strlem} we obtain $c\|\frac{x+y}{2}\|_\Phi>(1-\delta_1)|I_{D_{\delta_1}}|_\Phi=(1-\delta_1)c$. That is to say, we have both $\|x-y\|_\Phi\geq \lambda $ and $\|\frac{x+y}{2}\|_\Phi>(1-\delta_1)$, which is impossible since $x, y\in L^\Phi(E)$ and $\|x\|_\Phi=\|y\|_\Phi=1$.
\hfill\done

Notice that the assumption ``$(E,\|\cdot\|)$ is complete and has full support '' is used only in the process to show that ``$(L^\Phi(E),\|\cdot\|_\Phi)$ is uniformly convex'' implies that ``$(L^\Phi, |\cdot|_\Phi)$ is uniformly convex'', thus by removing this assumption and adding the assumption of uniform convexity of $(L^\Phi, |\cdot|_\Phi)$ we can obtain the following:

\begin{corollary}
 Let $(E,\|\cdot\|)$ be a complete RN module with base $(\Omega,{\mathcal F},P)$ and $\Phi$ a given Young function. If both $(L^\Phi, |\cdot|_\Phi)$ and $(L^\Phi(E),\|\cdot\|_\Phi)$ are uniformly convex, then $(E,\|\cdot\|)$ is random uniformly convex.
\end{corollary}

Since for every $p\in (1,\infty)$, $(L^p, |\cdot|_p)$ is uniformly convex, thus we have the following:

\begin{corollary}\label{puni}
Let $(E,\|\cdot\|)$ be an RN module and $1<p<\infty$. If $(L^p(E), \|\cdot\|_p)$ is uniformly convex then $(E,\|\cdot\|)$ is random uniformly convex.
\end{corollary}

This corollary is exactly Theorem 4.3 in  Guo and Zeng \cite{GZ1}.

The remainder is the most difficult part in this section. The technique in the proof of Proposition \ref{unisuf} follows the proof of Theorem 1(a) in Hudzik and Landes \cite{HD}. For our purpose, we first modify Lemma 1 of \cite{HD} as follows.
\begin{lemma}\label{discom}
Let $x,y$ be any two elements in an RN module $(E,\|\cdot\|)$ such that $P(B_{x,y})>0$, denote ${\hat x}=I_{[\|x\|\neq 0]}\frac{x}{\|x\|}, {\hat y}=I_{[\|y\|\neq 0]}\frac{y}{\|y\|}$, then the following statements hold:\\
\noindent(1). On the set $[~\|x\|>0~] \cap [~\|y\|>0~]$, we have that
$$  \|x+y\|  \leq  \big|\|x\|-\|y\|\big|+(\|x\|\wedge\|y\|)~\|{\hat x}+{\hat y}\| $$

\noindent(2). For any fixed real number $\gamma\in [0,2]$, on the set $[~\|x\|>0~] \cap [~\|y\|>0~]\cap [~\|{\hat x}+{\hat y}\|\leq \gamma~]$, we have that
$$\|x+y\|\leq \frac{\gamma}{2}(\|x\|+\|y\|)+(1-\frac{\gamma}{2})\big|\|x\|-\|y\|\big|$$
\noindent(3). If $(E,\|\cdot\|)$ is random uniformly convex, then for any real number $\epsilon\in (0,2]$, there exists a real number $\delta(\epsilon)\in (0,1)$ which is only decided by $\epsilon$ and independent of $x,y$ such that: for any real number $\eta\in (0,1]$, on the set $[~\|x\|>0~] \cap [~\|y\|>0~]\cap [~\|{\hat x}-{\hat y}\|\geq \epsilon~]\cap [~\|x\|\wedge\|y\| \geq \eta\{\|x\|\vee\|y\|\}~]$, we have
$$\|x+y\|\leq (1-\eta \delta(\epsilon))(\|x\|+\|y\|)$$
\end{lemma}

{\em proof.}

\noindent(1). On the set $[~\|x\|>0~] \cap [~\|y\|>0~]\cap [~\|x\|\geq \|y\|~]$, we have
\begin{eqnarray*}
\|x+y\|&= &\big\|(\|x\|-\|y\|){\hat x}+\|y\|({\hat x}+{\hat y})\big\| \\
       &\leq & \big\|(\|x\|-\|y\|){\hat x}\big\|+\big\|~\|y\|({\hat x}+{\hat y})~\big\|\\
       &=&\big|\|x\|-\|y\|\big| +\|y\|~\|{\hat x}+{\hat y}\| \\
       &=& \big|\|x\|-\|y\|\big|+(\|x\|\wedge\|y\|)~\|{\hat x}+{\hat y}\|,
\end{eqnarray*}
similarly, on the set $[~\|x\|>0~] \cap [~\|y\|>0~]\cap [~\|y\|\geq \|x\|~]$, we also have
$$\|x+y\|\leq \big|\|x\|-\|y\|\big|+(\|x\|\wedge\|y\|)~\|{\hat x}+{\hat y}\|$$
then (1) follows immediately.\\
\noindent(2). By using (1), on the set $[~\|x\|>0~] \cap [~\|y\|>0~]\cap [~\|{\hat x}+{\hat y}\|\leq \gamma~]$, we have
\begin{eqnarray*}
\|x+y\| & \leq & \big|\|x\|-\|y\|\big|+(\|x\|\wedge\|y\|)~\|{\hat x}+{\hat y}\| \\
        &\leq & \big|\|x\|-\|y\|\big|+\gamma~(\|x\|\wedge\|y\|~)\\
        &=& \big|\|x\|-\|y\|\big|+\gamma~ \{\frac{1}{2}(\|x\|+\|y\|)-\frac{1}{2}\big|\|x\|-\|y\|\big|\}\\
        &=& \frac{\gamma}{2}(\|x\|+\|y\|)+(1-\frac{\gamma}{2})\big|\|x\|-\|y\|\big|
\end{eqnarray*}
\noindent(3). By the random uniformly convexity of $(E,\|\cdot\|)$, for any real number $\epsilon\in (0,2]$, since we can take $\epsilon$ as an element in $L_\#[0,2]$, there exists a $\lambda(\epsilon)\in L_\#[0,1]$ which is only decided by $\epsilon$ such that: on the set $[~\|x\|>0~] \cap [~\|y\|>0~]\cap [~\|{\hat x}-{\hat y}\|\geq \epsilon~]$ we have $\|{\hat x}+{\hat y}\|\leq 2-2\lambda(\epsilon)$, by the meaning of $L_\#[0,1]$, we can choose a real number $\delta(\epsilon)\in (0,1)$ such that $\delta(\epsilon)\leq\lambda(\epsilon)$. Therefore it follows from (1) that, on the set $[~\|x\|>0~] \cap [~\|y\|>0~]\cap [~\|{\hat x}-{\hat y}\|\geq \epsilon~]\cap [~\|x\|\wedge\|y\| \geq \eta\{\|x\|\vee\|y\|\}~]$, we have
\begin{eqnarray*}
  \|x+y\| & \leqslant & \big|\|x\|-\|y\|\big|+(\|x\|\wedge\|y\|)~\|{\hat x}+{\hat y}\| \\
          &\leq & \big|\|x\|-\|y\|\big|+(\|x\|\wedge\|y\|)(2-2\delta(\epsilon)) \\
          &=& \|x\|+\|y\| -2\lambda(\epsilon)(\|x\|\wedge\|y\|)\\
          &\leq & \|x\|+\|y\| -2\lambda(\epsilon)\eta (\|x\|\vee\|y\|)\\
          &\leq & \|x\|+\|y\| -\eta\lambda(\epsilon)(\|x\|+\|y\|)~~(~\text{since}~ 2(\|x\|\vee\|y\|)\geq \|x\|+\|y\|~)\\
          &=&(1-\eta\lambda(\epsilon))(\|x\|+\|y\|) \\
          &\leq &(1-\eta\delta(\epsilon))(\|x\|+\|y\|)
\end{eqnarray*}
\hfill\done

Now we can state and prove Proposition \ref{unisuf}. We remind that the norm $|\cdot|_\Phi$ on the Orlicz space $(L^\Phi, |\cdot|_\Phi)$ is a lattice norm, namely (i) for any $\xi\in L^\Phi$, $|\xi|_\Phi=\big|~|\xi|~\big|_\Phi$,  and (ii) for any $\xi,\eta\in L^\Phi$ with $0\leq \xi\leq \eta$, $|\xi|_\Phi\leq |\eta|_\Phi$.

\begin{proposition}\label{unisuf}
  Let $(E,\|\cdot\|)$ be a complete RN module with base $(\Omega,{\mathcal F},P)$ and $\Phi$ a given Young function. If $(L^\Phi, |\cdot|_\Phi)$ is uniformly convex and $(E,\|\cdot\|)$ is random uniformly convex, then $(L^\Phi(E), \|\cdot\|_\Phi)$ is uniformly convex.
\end{proposition}

{\em proof.}
We use an equivalent definition of uniform convexity. For any two sequences $\{x_n, n\in {\mathbb N}\}$ and $\{y_n, n\in {\mathbb N}\}$ in $(L^\Phi(E), \|\cdot\|_\Phi)$ with $\|x_n\|_\Phi=\|y_n\|_\Phi=1, \forall n\in {\mathbb N}$ such that $\|x_n+y_n\|_\Phi\to 2$ as $n\to\infty$, we must show $\|x_n-y_n\|_\Phi\to 0$ as $n\to\infty$.

For each $n\in {\mathbb N}$, let $u_n=\|x_n\|, v_n=\|y_n\|, s_n=\|x_n+y_n\|, S_n=u_n+v_n, d_n=\|x_n-y_n\|$ and $D_n=|u_n-v_n|$. Since $ 2=|u_n|_\Phi+|v_n|_\Phi\geq |S_n|_\Phi\geq |s_n|_\Phi$ and $|s_n|_\Phi\to 2$ as $n\to\infty$, we have that $|S_n|_\Phi-|s_n|_\Phi\to 0$ as $n\to\infty$.

Choose a sequence of positive real numbers $\{\epsilon_k,k\in {\mathbb N}\} $ such that $\epsilon_k\to 0$ as $k\to\infty$. For each $k\in {\mathbb N}$, let $\delta(\epsilon_k)$ as (3) of Lemma \ref{discom}, then there exists a subsequence of positive integers $N_1<N_2< \cdots <N_k<\cdots $ such that $|S_n|_\Phi-|s_n|_\Phi\leq \frac{1}{k}\delta(\epsilon_k), \forall n\geq N_k$. Then we can take two sequences of positive numbers $\{\eta_n, n\in {\mathbb N}\}$ and $\{\varepsilon_n, n\in {\mathbb N}\}$ as follows: for $n<N_1$, let $\eta_n=1$ and $\varepsilon_n=1$; for $N_k\leq n<N_{k+1}$, choose $\eta_n=\frac{1}{\sqrt k}$ and $\varepsilon_n=\epsilon_k, k=1,2,\dots$. Obviously we have $\eta_n\to 0,~\varepsilon_n\to 0$ and $(|S_n|_\Phi-|s_n|_\Phi)/\gamma_n\to 0$ as $n\to\infty$, where $\gamma_n=\eta_n \delta(\varepsilon_n)$ for each $n\in {\mathbb N}$.

(a). Define the sets $A(\eta_n)=[~u_n>0~]\cap [~v_n>0~]\cap [~(u_n\wedge v_n )\geq \eta_n(u_n\vee v_n)~]$ and $A_n=A(\eta_n)\cap [~\|{\hat x_n}-{\hat y_n}\|\geq \varepsilon_n]$, $B_n=A(\eta_n)\setminus A_n$ for each $n\in {\mathbb N}$.

From (3) of Lemma \ref{discom}, we have that on the set $A_n$, $s_n\leq (1-\gamma_n)S_n$. Then on $\Omega$, $s_n\leq S_n(1-\gamma_nI_{A_n})$, therefore $|s_n|_\Phi\leq |S_n(1-\gamma_nI_{A_n})|_\Phi=|S_n(1-\gamma_n+\gamma_nI_{A_n^c})|_\Phi\leq (1-\gamma_n)|S_n|_\Phi+\gamma_n|S_nI_{A_n^c}|_\Phi$,
which yields that $|S_nI_{A_n^c}|_\Phi\geq |S_n|_\Phi-(|S_n|_\Phi-|s_n|_\Phi)/{\gamma_n}\to 2$ as $n\to \infty$, clearly $|S_nI_{A_n^c}|_\Phi\leq |S_n|_\Phi\leq 2$, thus we conclude that $|S_nI_{A_n^c}|_\Phi\to 2$ as $n\to \infty$.

(b). Applying (2) of Lemma \ref{discom}  to $x_n$ and $-y_n$, we conclude that on $B_n$, $d_n\leq \frac{\varepsilon_n}{2} S_n+(1-\frac{\varepsilon_n}{2})D_n$, so that on $\Omega$,
\begin{eqnarray*}
 d_n &\leq & d_nI_{B_n^c}+[\frac{\varepsilon_n}{2} S_n+(1-\frac{\varepsilon_n}{2})D_n]I_{B_n}\\
      &=& \frac{\varepsilon_n}{2}[d_nI_{B_n^c}+S_nI_{B_n}]+(1-\frac{\varepsilon_n}{2})[d_nI_{B_n^c}+D_nI_{B_n}]\\
      &\leq & \frac{\varepsilon_n}{2}S_n +(1-\frac{\varepsilon_n}{2})[d_nI_{B_n^c}+D_nI_{B_n}],
\end{eqnarray*}
thus $|d_n|_\Phi\leq \frac{\varepsilon_n}{2}|S_n|_\Phi+(1-\frac{\varepsilon_n}{2})|d_nI_{B_n^c}+D_nI_{B_n}|_\Phi\leq \varepsilon_n +(1-\frac{\varepsilon_n}{2})|d_nI_{B_n^c}+D_nI_{B_n}|_\Phi$.
Since $d_nI_{B_n^c}+D_nI_{B_n}\leq S_nI_{B_n^c}+D_nI_{B_n}$, we obtain that $$|S_nI_{B_n^c}+D_nI_{B_n}|_\Phi\geq |d_nI_{B_n^c}+D_nI_{B_n}|_\Phi\geq \frac{|d_n|_\Phi-\varepsilon_n}{1-\varepsilon_n/2}.$$

(c). On $A^c(\eta_n)=\Omega\setminus A(\eta_n)$, $S_n-D_n=u_n+v_n-|u_n-v_n|=2(u_n\wedge v_n)\leq 2\eta_n(u_n+v_n)$, thus  $|(S_n-D_n)I_{A^c(\eta_n)}|_\Phi\leq 2\eta_n(|u_n|_\Phi+|v_n|_\Phi)= 4\eta_n\to 0$ as $n\to \infty$.

(d). Define $\xi_n=v_n(I_{A_n^c}-I_{A_n})$, then $|\xi_n|=v_n$ so that $|\xi_n|_\Phi=|v_n|_\Phi=1$. Since $|u_n+\xi_n|=|(u_n+v_n)I_{A_n^c}+(u_n-v_n)I_{A_n}|\geq (u_n+v_n)I_{A_n^c}=S_nI_{A_n^c}$, thus $|u_n+\xi_n|_\Phi\geq |S_nI_{A_n^c}|_\Phi\to 2$ as $n\to \infty$, then by the uniform convexity of $(L^\Phi, |\cdot|_\Phi)$, we must have $|u_n-\xi_n|_\Phi\to 0$ as $n\to \infty$. However
\begin{eqnarray*}
 |u_n-\xi_n| & = & |(u_n+v_n)I_{A_n}+(u_n-v_n)I_{A_n^c}| \\
           & = & |(u_n+v_n)|I_{A_n}+|(u_n-v_n)|I_{A_n^c} \\
           & = & S_n I_{A_n}+D_n I_{A_n^c} \\
           & = & S_n(I_{B_n^c}-I_{A^c(\eta_n)})+D_n(I_{B_n}+I_{A^c(\eta_n)})\\
           & = & S_n I_{B_n^c}+D_n I_{B_n}- (S_n-D_n)I_{A^c(\eta_n)},
\end{eqnarray*}
thus
$|u_n-\xi_n|_\Phi\geq |S_nI_{B_n^c}+D_nI_{B_n}|_\Phi-|(S_n-D_n)I_{A^c(\eta_n)}|_\Phi$,

therefore,
\begin{eqnarray*}
\limsup_{n\to\infty} |d_n|_\Phi &=& \limsup_{n\to\infty} \frac{|d_n|_\Phi-\varepsilon_n}{1-\varepsilon_n/2}~~(~\text{since}~\lim_{n\to\infty}\varepsilon_n=0~)\\
                     &\leq & \limsup_{n\to\infty} |S_nI_{B_n^c}+D_nI_{B_n}|_\Phi ~~(~\text{from part (b)})\\
                     &\leq & \limsup_{n\to\infty}~(|u_n-\xi_n|_\Phi+|(S_n-D_n)I_{A^c(\eta_n)}|_\Phi)~~(~\text{from part (d)}) \\
                     &=& 0+0=0 ~~(~\text{from part (c) and part (d)})
\end{eqnarray*}

\hfill\done

Since for every $p\in (1,\infty)$, $(L^p, |\cdot|_p)$ is uniformly convex, thus we have the following:

\begin{corollary}
Let $(E,\|\cdot\|)$ be an RN module and $1<p<\infty$. If $(E,\|\cdot\|)$ is random uniformly convex, then $(L^p(E), \|\cdot\|_p)$ is uniformly convex.
\end{corollary}

This combines Corollary \ref{puni} is exactly Theorem 1.8 of Guo and Zeng \cite{GZ2}.

 \vspace{0.8cm}
{\noindent\bf Acknowledgements.}
The work of the first author was supported by the Natural Science Foundation of China (No.11701531) and the Fundamental Research Funds for the Central Universities, China University of Geosciences (Wuhan). The work of the second author was supported by
the Natural Science Foundation of China (Grant No. 11501580). The work of the third author was supported by the Natural Science Foundation of China (Grant No. 11301568) and the Chongqing Technology and Business University Scientific Research Foundation 2012-56-10.


\end{document}